\documentstyle[amsfonts,aps]{revtex}

\begin{document}
\title{NON-COMMUTATIVE COREPRESENTATIONS OF QUANTUM GROUPS\\
\smallskip\ }
\author{H. Montani$^{*}${\bf \ }}
\address{\smallskip\ \ \\
{\it Instituto Balseiro - Centro At\'{o}mico Bariloche}\\
{\it 8400 - S. C. de Bariloche, Rio Negro, Argentina.}\\
\medskip}
\author{R. Trinchero$^{\dagger }$}
\address{\smallskip\ \\
{\it Centre de Physique Th\'{e}orique - CNRS - Luminy, Case \thinspace
\thinspace \thinspace 907, }\\
{\it F-13288 Marseille Cedex 9 - France.}\\
{\it \ }}
\maketitle

\begin{abstract}
We consider a twisted version of quantum groups corepresentations. This
generalization amounts to include in the theory the case where quantum space
coordinates and its endomorphism matrix entries belong to a non-commutative
quadratic algebra.
\end{abstract}

\vspace*{\fill}

$*$e-mail: montani@cab.cnea.edu.ar

$\dagger $On leave of absence from {\it Instituto Balseiro - Centro
At\'{o}mico Bariloche, 8400- S. C. de Bariloche, Rio Negro, Argentina. }

e-mail: trincher@cab.cnea.edu.ar
\vskip 0.1cm
\hrule\hfill CPT-97/P.3498

\newpage\

\section{Introduction}

Quantum Groups arise as the abstract structure underlying the symmetries of
integrable systems in (1+1) dimensions\cite{fad}. There, the theory of
quantum inverse scattering give rise to some deformed algebraic structures
which were first explained by Drinfel'd as deformations of classical Lie
algebras\cite{drinf}\cite{jim}. An analog structure was obtained by
Woronowicz in the context of non-commutative $C^{*}$-algebras\cite{wor}.
There is a third approach, due to Manin, where Quantum Groups are
interpreted as the endomorphisms of certain non commutative algebraic
varieties defined by quadratic algebras, called quantum linear spaces (QLS)%
\cite{man}. Faddeev {\it et al }had also interpreted the Quantum Groups from
the point of view of corepresentations and quantum spaces, furnishing a
connection with the quantum deformations of the universal enveloping
algebras and the quantum double of Hopf algebras \cite{frt}\cite{kass}.

From the algebraic point of view, quantum groups are Hopf algebras and the
relation with the endomorphism algebra of QLS come from their
corepresentations on tensor product spaces. The usual construction of the
coaction on the tensor product space involves the flip operator
interchanging factors of the tensor product of the QLS with the bialgebra.
This fact implies the commutativity between the matrix elements of a
representation of the endomorphism and the coordinates of the QLS. Moreover,
the flip operator for the tensor product is also involved in many steps of
the construction of Quantum Groups. In the braided approach to $q$%
-deformations, the flip operator is replaced by a braiding giving rise to
the quasi-tensor category of $k$-modules, where a natural braided coaction
appears\cite{maj}.

In the present work, we introduce a twisted coaction over the tensor product
space, thus admitting non-commutative relations between endomorphism matrix
entries and quantum linear space coordinates, however this has nothing to do
with the braided approach mentioned above. We find the conditions under
which the general algebraic framework of multiplicative quantum groups still
holds. It is also shown that the bialgebras arising from this context may be
regarded as a partial twisting of usual quantum groups and the connections
with integrable systems is analyzed after the introduction of the spectral
parameter. This twisted coaction allow us to introduce new deformation
parameters in the endomorphism bialgebra of the QLS, as it is shown in the
quantum plane example where a four parameters deformation is obtained,
although the Yang-Baxter condition is relaxed. Also, we find a non central
object playing the role of $(q,p,r,s)$-deformed determinant. In the
undeformed limit for the parameters $(r,s)$ we recover the biparametric
deformation $GL_{qp}(2)$ described in ref. \cite{man2}.

We present a brief description of the corepresentations of bialgebras in the
second section, developing our approach to modified corepresentations in the
third one. In the fourth section we present the result of the previous
section as twisted bialgebras. The connection with integrable systems is
discussed in the fifth section and, finally, we work out the quantum plane
example in the last section.

\section{Quantum algebras and corepresentations}

Let $V$ be a vector space of dimension $n$, $\{e_{i}\}$ a basis for $V$ and $%
H_{o}$ the trivial bialgebra of functions over $GL(n,{\bf C})$. This
bialgebra is freely generated by the identity and the coordinates functions $%
T_{i}^{j}$, in the basis $\{e_{i}\}$, defined by

\[
T{_i^j:}GL(n,C)\longrightarrow C
\]

\[
T{_i^j:\,}g\longrightarrow g_i^j
\]
for $g\in GL(n,{\bf C})$. The $T{_i^j}$ are group like, the coproduct and
counit are given by

\begin{equation}
\Delta T_i^j=T_i^k\otimes T_k^j  \label{comult}
\end{equation}

\begin{equation}
\varepsilon (T_i^j)=\delta _j^i  \label{coun}
\end{equation}
From now on, summation over repeated index is assumed. The comodule $(\delta
,V)$, with

\[
{\delta :}V{\longrightarrow }H_o{{\otimes }V}
\]

\begin{equation}
\delta (e_{i})=T_{i}^{j}\otimes e_{j}  \label{coac}
\end{equation}
provides a representation of $GL(n,{\bf C})$ in $V$, through the $g_{i}^{j}$
in the basis $b$ of $V$. It has the coassociativity property and preserves
the counit, which is expressed by the relations

\begin{equation}
(I_{H_{o}}\otimes \delta _{V})\,\delta _{V}=\left( \Delta \otimes
I_{V}\right) \,  \label{coass}
\end{equation}

\begin{equation}
\left( \epsilon \otimes I_{H_{o}}\right) \delta _{V}=I_{V}  \label{counit}
\end{equation}
In order to extend the comodule to the tensor product algebra $V^{\otimes },$
one introduce the coaction on $V\otimes V,$

\begin{equation}
\begin{array}{c}
\delta _{V\otimes V}:V\otimes V\longrightarrow H_o\otimes V\otimes V \\
\\
\delta _{V\otimes V}=(m\otimes I_{V\otimes V})\,(I_{H_o}\otimes \tau \otimes
I_V)\,(\delta _V\otimes \delta _V)
\end{array}
\label{coacvv}
\end{equation}
The extension to $V^{\otimes N}$ is achieved via the recursive relations

\begin{equation}
\begin{array}{c}
\delta _{V^{\otimes N}}:V^{\otimes N}\longrightarrow H_{o}\otimes V^{\otimes
N} \\
\\
\delta _{V^{\otimes N}}=(m\otimes I_{V^{\otimes N}})\,(I_{H_{o}}\otimes \tau
_{V^{\otimes N},H_{o}}\otimes I_{V})\,(\delta _{V^{\otimes (N-1)}}\otimes
\delta _{V})
\end{array}
\label{vncoact}
\end{equation}
where $\tau _{V^{\otimes N},H_{o}}$ is the flip operator mapping $\tau
_{V^{\otimes N},H_{o}}:V^{\otimes N}\otimes H_{o}\longrightarrow
H_{o}\otimes V^{\otimes N}.$ This definition satisfy the coassociativity and
counit properties described in relations (\ref{coass}) and (\ref{counit}).
It is also worth to remark that the appearance of the flip operator in (\ref
{vncoact}) leads to the commutativity between the coordinates of the quantum
space and its endomorphism matrix entries , as its assumed for the quantum
plane and $GL_{q}(2)$ \cite{man}.

Building up corepresentations for objects with more structure than $%
V^{\otimes }$, as quadratic algebras for example, requires some extra
conditions that we sketch below.

Let $A$ denote the quadratic algebra generated by the ideal $I({\frak B}),$
where ${\frak B}:V\otimes V\longrightarrow V\otimes V$, then

\begin{equation}
A({\frak B})=\frac{V^{\otimes }}{I({\frak B})}  \label{alq}
\end{equation}
$V^{\otimes }$ is tensor algebra on $V.$ In general we consider ${\frak B}$
with the form
\begin{equation}
{\frak B}=(I_{V\otimes V}-B)  \label{bb}
\end{equation}

\begin{equation}
e_{i}e_{j}-B_{ij}^{kl}e_{k}e_{l}  \label{eebb}
\end{equation}
$\delta _{V\otimes V}$ must be an homomorfism of quadratic algebra, i.e.,

\begin{equation}
\left( I_{H_{o}}\otimes {\frak B}\right) \delta _{V\otimes V}=\delta
_{V\otimes V}\,{\frak B}  \label{bddb}
\end{equation}

This is satisfied if $H$ is the bialgebra arising from the quotient of the
free algebra generated by the objects $T_{i}^{j}$ and the ideal $I({\frak B}%
,H_{o})$ generated by the quadratic relation

\begin{equation}
B_{ij}^{kl}T_k^rT_l^s-T_i^kT_j^lB_{kl}^{rs}  \label{btt}
\end{equation}
i.e.,

\begin{equation}
H=\frac{H_o}{I({\frak B},H_o)}  \label{hbtt}
\end{equation}

Since $I({\frak B},H_o)$ is a coideal with relation to $\Delta $ , $H$
becomes a bialgebra, namely an {\it FRT bialgebra}. The equation (\ref{btt})
is a central object in the so called FRT construction\cite{frt}. In this
way, $A({\frak B})$ becomes in a $H$-algebra comodule.

\section{Generalized corepresentations}

The main aim of this section is to build up the mathematical framework
encoding the situation in which entries of the endomorphism matrix may not
commute with the coordinates of the quantum linear space defined in (\ref
{alq}). We will reach it by means of a modification in the corepresentation
theory, obtained by substituting the flip map $\tau $ in the standard
definition of the coaction on $V\otimes V$ by a non trivial map $\gamma .$

As described in the previous section, supplying the quantum linear space
with a comodule structure requires a right definition of a coaction on the
tensor product space, and the standard definition of $\delta _{V\otimes V}$,
eq. (\ref{coacvv}), provides both $V\otimes V$ with a $H_o$-comodule
structure and $A=\frac{V^{\otimes }}{I({\frak B})}$ with an $H$-comodule
structure.

Then, let us introduce the map $\gamma $ , defined by

\begin{equation}
\gamma :V\otimes H\longrightarrow H\otimes V  \label{g}
\end{equation}

\begin{equation}
\gamma (e_i\otimes T{_j^k)=\gamma }_{ijn}^{klm}T_l^n\otimes e_m  \label{egt}
\end{equation}
and our proposal of generalized or non-commutative coaction on tensor
product space is
\begin{equation}
\delta _{V\otimes V}^\gamma =(m\otimes I_{V\otimes V})\,(I_{H_\gamma
}\otimes \gamma \otimes I_V)\,(\delta _V\otimes \delta _V)  \label{coacg}
\end{equation}
Then, we shall see that a $H_\gamma $-comodule structure there it is
possible, for some $H_\gamma $ to be constructed and provided $\gamma $
satisfying some requirements. The first question is finding the condition
under which $\delta _{V\otimes V}^\gamma $ is actually a coaction.\ It is
addressed in the following proposition:
\vskip 0.3cm
\noindent{\bf Proposition 1:} {\it The map $\delta _{V\otimes V}^{\gamma
}:V\otimes V\rightarrow H_{o}\otimes V\otimes V$ is a coaction turning
$V\otimes V$ into a $H_{o}$-comodule iff $ \gamma :V\otimes
H_{o}\longrightarrow H_{o}\otimes V$ satisfies the following conditions}
\begin{equation}
{\gamma }_{ijn}^{klm}=\delta _{i}^{m}\theta _{jn}^{kl}  \label{gamarel1}
\end{equation}

\begin{equation}
\theta _{ij}^{pl}\theta _{pk}^{rs}-\delta _{j}^{s}\theta _{ik}^{rl}=0
\label{gamarel2}
\end{equation}

\begin{equation}
\theta _{jn}^{kn}=\delta _{j}^{k}  \label{gamarel3}
\end{equation}
.

{\it Proof: } These properties for $\gamma $ arise straightforward from the
coassociativity and counit conditions

\begin{equation}
\begin{array}{c}
(\Delta \otimes I_{V\otimes V})\circ \delta _{V\otimes V}^{\gamma
}=(I_{H_{\gamma }}\otimes \delta _{V\otimes V}^{\gamma })\circ \delta
_{V\otimes V}^{\gamma } \\
\\
(\epsilon \otimes I_{V\otimes V})\circ \delta _{V\otimes V}^{\gamma
}=I_{V\otimes V}
\end{array}
\label{coaction}
\end{equation}
$\Box $
\vskip 0.3cm

A mapping ${\gamma }$ satisfying the conditions (\ref{gamarel1}-\ref
{gamarel3}) leads to a comodule over $V^{\otimes N}$ as stated in the
following proposition.
\vskip 0.3cm

\noindent{\bf Proposition 2:} {\it Let $\delta _{V^{\otimes N}}:V^{\otimes
N} \longrightarrow H_{o}\otimes V^{\otimes N}$ be defined by,
\begin{equation}
\delta _{V^{\otimes N}}^{\gamma }=(m_{H_{o}}\otimes I_{V^{\otimes
N}})\,(I_{H_{o}}\otimes \gamma _{V^{\otimes (N-1)},H_{o}}\otimes
I_{V})\,(\delta _{V^{\otimes (N-1)}}\otimes \delta _{V})  \label{gamacoac}
\end{equation}
where
\[
\gamma _{V^{\otimes N},H_{o}}(e_{i_{1}}\otimes ...\otimes e_{i_{N}}\otimes
T_{j}^{k})=\gamma _{V,H_{o}}\otimes I_{V^{\otimes (N-1)}}\,(e_{i_{1}}\otimes
\gamma _{V^{\otimes (N-1)},H_{o}}\,(e_{i_{2}}\otimes ...\otimes
e_{i_{N}}\otimes T_{j}^{k}))
\]
and
\[
\gamma _{V,H_{o}}=\gamma
\]
then $(V^{\otimes }\,,$ $\{\delta _{V^{\otimes N}}^{\gamma }\})$ is a left $%
H_{o}$-comodule.}
\vskip 0.3cm

The proof runs as the previous one, just with a more complicated algebra.

Recalling the bijection between comodules and multiplicative matrices\cite
{man}, let us consider the multiplicative matrix $M$ in $V\otimes V$ with
coefficients in $H_o$ corresponding to the comodule $\delta _{V\otimes
V}^\gamma $ ,i.e.,

\begin{equation}
\begin{array}{c}
\delta _{V\otimes V}^{\gamma }\equiv M\in End(V\otimes V,H_{o}) \\
\\
\delta _{V\otimes V}^{\gamma }(e_{i}\otimes e_{j})=M_{ij}^{rs}\otimes
e_{r}\otimes e_{s}
\end{array}
\label{M0}
\end{equation}
hence $M$ is

\begin{equation}
M_{ij}^{kl}=T_{i}^{k}\theta _{jn}^{lm}T_{m}^{n}.  \label{M}
\end{equation}
Let us adopt the following convention: for $A_{ij}^{kl}$ and $D_{ij}^{kl}$
being any pair of four-tensors, we write $(A\times
B)_{ij}^{rs}=A_{ij}^{kl}\times D_{kl}^{rs}$, where $\times $ stands for any
kind of product (tensor, algebraic,etc.), and sum over repeated index is
also assumed.

With this notation,conditions (\ref{gamarel2}) and (\ref{gamarel3}) are

\begin{equation}
\begin{array}{c}
\Delta M=M\otimes M \\
\\
\epsilon \,(M)=I
\end{array}
\label{dmmm}
\end{equation}
The next step is to consider a quadratic structure on $V\otimes V$ giving
rise to a QLS$.$ Now, the bialgebra $H_o$ is no longer in the endomorphism
algebra of the QLS. Let us consider a QLS generated by the quotient algebra

\begin{equation}
A({\frak B})=\frac{V^{\otimes }}{I({\frak B})},  \label{alq2}
\end{equation}
where ${\frak B}$ means the relations defining the quadratic algebra.
Associated with it we now introduce a new bialgebra structure on the free
algebra generated by the $\{T_{i}^{k}\}.$
\vskip 0.3cm

\noindent{\bf Proposition 3:} {\it
Let $H_{o}$ the free algebra generated by the $\{T_{i}^{k}\}$, $\gamma $ as
in the previous proposition and $I({\frak B}M-M{\frak B})$ is the ideal
generated by the quadratic relation
\begin{equation}
{\frak B}M-M{\frak B}  \label{bmb}
\end{equation}
then, the quotient algebra $H_{\gamma }$ defined as
\begin{equation}
H_{\gamma }=\frac{H_{o}}{I({\frak B\,}M-M\,{\frak B})}  \label{hgam}
\end{equation}
is a bialgebra.}

{\it Proof: }A necessary and sufficient condition for $H_\gamma $ to be a
bialgebra is that $I({\frak B\,}M-M\,{\frak B})$ be a coideal,i.e.,
\[
\Delta I\subset I\otimes H_o+H_o\otimes I
\]
Then, taking into account the relation (\ref{dmmm}), one gets

\begin{eqnarray*}
\Delta {(}{\frak B}M-M{\frak B}) &=&{\frak B\,}\Delta M-\Delta M\,{\frak B}
\\
&=&{\frak B\,}(M\otimes M)-(M\otimes M)\,{\frak B} \\
&=&{(}{\frak B\,}M-M\,{\frak B})\otimes M-M\otimes {(}{\frak B\,}M-M\,{\frak %
B})
\end{eqnarray*}
and

\[
\epsilon {(}{\frak B\,}M-M\,{\frak B})={\frak B\,}\epsilon (M)-\epsilon (M)\,%
{\frak B}=0
\]
hence $H_\gamma $ is a bialgebra.$\Box $
\vskip 0.3cm

The main result of this section is expressed in the following proposition
\vskip 0.3cm

\noindent{\bf Proposition 4:} {\it
$\{\delta _{V^{\otimes N}}^{\gamma }\}$ supplies $A({\frak B})=%
{\displaystyle {V^{\otimes } \over I({\frak B})}}%
$ with an left $H_{\gamma }$-comodule structure.}

{\it Proof:} This assertion means the map $\delta _{V\otimes V}^{\gamma }$ :$%
A({\frak B})\rightarrow H_{\gamma }\otimes A({\frak B})$ is an homomorphism
of quadratic algebras, as in eq. (\ref{bddb}). This fact is realized by the
commutation relation

\[
(I_{H_\gamma }\otimes {\frak B})\circ \delta _{V\otimes V}^\gamma =\delta
_{V\otimes V}^\gamma \circ {\frak B}
\]
which is immediately satisfied by virtue of the ideal defining $H_\gamma ,$
i.e., the condition

\[
{\frak B\,}M=M\,{\frak B}
\]

$\Box $
\vskip 0.3cm

Resuming, we can make the following assertion: given a quadratic algebra $A(%
{\frak B})$ and a map $\gamma $ satisfying the relations (\ref{gamarel1}-\ref
{gamarel3}), then $H_{\gamma }=\frac{H_{o}}{I({\frak B}M-M{\frak B})}$ is a
bialgebra and $\{\delta _{V^{\otimes N}}^{\gamma }\}$ renders $A({\frak B})$
into a $H_{\gamma }$-comodule (a similar quadratic algebra arise in the
context of quantum braided group\cite{hlavaty}). This may be understood
because of the bijection between all the structures of left comodule on $V=%
{\bf C}^{n}$ and the multiplicative matrix ${\bf M(}n,H_{\gamma })$\cite{man}%
, since $M\in H_{\gamma }$ satisfy $\Delta M=M\otimes M$ and $\varepsilon
(M)=I_{V\otimes V},$ for $M\in H_{\gamma }$ .

The existence of an antipode is not involved in the comodule structure, so
we may expect the above construction still holds when $H_{\gamma }$ is a
Hopf algebra, giving rise to non-commutative corepresentations of Quantum
Groups.

In the last section, we describe an explicit example enjoying all these
properties presented above, namely a multiparameter deformed version of the
endomorphism of the quantum plane.

\section{Relations with twisted bialgebras}

Let us introduce a bialgebra structure on $Hom(H_{o}^{\otimes 2},k)$ by
means the convolution product $*$ of linear forms, defined as $%
(f*g)(T)=(f\otimes g)(\Delta T)$ for $f,g$ $\in Hom(H_{o}^{\otimes 2},k)$
and $T$ $\in H_{o}^{\otimes 2}.$ The coproduct is $(\Delta h)(T\otimes
T^{\prime })=(h\circ m)(T\otimes T^{\prime })$ for $h\in Hom(H_{o}^{\otimes
2},k)$ and $T,T^{\prime }\in H_{o}.$ The unit is $\varepsilon ,$ the counit
of the $H_{o}$ , namely $(\varepsilon *$ $f)(T)=f(T).$

In this framework \cite{drinf}\cite{maj}\cite{kass}, $H$ defined in eq. (\ref
{hbtt}) can be presented as the bialgebra $H(m,\Delta ,\eta ,\varepsilon
,R), $ with $R:$ $H_o^{\otimes 2}\rightarrow k$ being an invertible linear
form, related to ${\frak B}$ of the previous section by

\begin{equation}
R(T_i^k\otimes T_j^l)=R_{ij}^{kl}\equiv B_{ji}^{kl}  \label{r}
\end{equation}
and defined by the quadratic ideal generated by the relation
\begin{equation}
m^{op}=R*m*\overline{R}  \label{rtt1}
\end{equation}
Here, $m^{op}=m\circ \tau $ and $R*\overline{R}=\overline{R}*R=\varepsilon .
$ Moreover, $H$ is said {\it dual quasi-triangular}\cite{maj} provided $R$
satisfies
\begin{equation}
\begin{array}{rrr}
R\circ (I_H\otimes m)=R_{13}*R_{12} &  & R\circ (m\otimes I_H)=R_{13}*R_{23}
\end{array}
\label{rm}
\end{equation}
Here, $R_{12}=R\otimes \varepsilon $ , $R_{23}=\varepsilon \otimes R$ , $%
R_{13}=(\varepsilon \otimes R)\circ (\tau \otimes I_H).$ This last relations
implies $R$ is a solution of the {\it Quantum Yang-Baxter equation}

\begin{equation}
R_{12}*R_{13}*R_{23}=R_{23}*R_{13}*R_{12}  \label{YB}
\end{equation}

Coming back to our problem, let us work out the bialgebra structure $%
H_\gamma $, eq.(\ref{hgam}), derived from non-commutative corepresentations
of the previous section. The following characterization of the $\gamma $ map
drives to a different interpretation of the bialgebra $H_\gamma .$

Let $\theta $ be the linear map
\begin{equation}
\begin{array}{c}
\theta :H_o\longrightarrow H_o \\
\\
\theta (T_i^j)=\theta _{im}^{jn}\,T_n^m=\widetilde{T}_i^j
\end{array}
\label{hmap}
\end{equation}
with the properties
\begin{equation}
\begin{array}{c}
\theta _{ij}^{pl}\theta _{pk}^{rs}-\delta _j^s\theta _{ik}^{rl}=0 \\
\\
\theta _{jn}^{kn}=\delta _j^k
\end{array}
\label{tita}
\end{equation}
\vfill\eject
\noindent{\bf Proposition 5:} {\it
 $\theta $ is a coalgebra homomorphism, such that}

\begin{equation}
\begin{array}{c}
\Delta \widetilde{T}_{i}^{j}=\widetilde{T}_{i}^{k}\otimes \widetilde{T}%
_{k}^{j} \\
\\
\epsilon \,(\widetilde{T}_{i}^{j})=\delta _{i}^{j}
\end{array}
\label{copts}
\end{equation}
{\it Proof: }the properties (\ref{tita}) implies that $\widetilde{T}$ are
group like elements, then satisfying the coassociativity and counit
properties.
\vskip 0.3cm

With this notation, the twisting $\gamma $, eq. (\ref{g}) can be
now expressed as \[
\gamma (e_i\otimes T_j^k)=\theta (T_j^k)\otimes e_i=\widetilde{T}_j^k\otimes
e_i
\]
and the quadratic relation (\ref{bmb}) can be written more explicitly as

\begin{equation}
R_{ij}^{ab}\,T_a^k\,\widetilde{T}_b^k=T_j^a\,\widetilde{T}_i^b\,R_{ba}^{kl}
\label{rtt}
\end{equation}
This relation generates the ideal which give rise to the quadratic algebra $%
H_\gamma (m,\Delta ,\eta ,\varepsilon ).$
\vskip 0.3cm

\noindent {\bf Proposition 7:} {\it
Let $\theta $ be an automorphism in $H_{o}$, then there is an isomorphism of
the bialgebra $H_{\gamma }(m,\Delta ,\eta ,\varepsilon )$ with the bialgebra
$H(m_{\theta },\Delta ,\eta ,\varepsilon ,R)$, where the deformed product $%
m_{\theta }:H_{o}\otimes H_{o}\rightarrow H_{o}$ is defined as
\begin{equation}
\begin{array}{c}
m_{\theta }=m\circ (I_{H}\otimes \theta ) \\
\\
(I_{H}\otimes \theta )\,(T_{i}^{k}\otimes T_{j}^{l})=T_{i}^{k}\,\otimes
\widetilde{T}_{j}^{l}
\end{array}
\,\,\,,  \label{mtita}
\end{equation}
and it fullfil the relation
\begin{equation}
m_{\theta }^{op}\,=R*m_{\theta }*\,\overline{R}  \label{rtt2}
\end{equation}}

{\it Proof: }The ideal generated by (\ref{rtt2}) is exactly (\ref{rtt}), and
since $\theta $ is an automorphism the proof is obvious.. Associativity can
be reached, for example, with
\begin{equation}
\begin{array}{c}
m_{\theta }(T_{i}^{l}\,T_{j}^{m}\otimes T_{k}^{n})=T_{i}^{l}\,T_{j}^{m}%
\widetilde{\widetilde{T}}_{k}^{n} \\
\\
m_{\theta }(T_{i}^{l}\,\otimes T_{j}^{m}T_{k}^{n})=T_{i}^{l}\,\widetilde{T}%
_{j}^{m}\widetilde{T}_{k}^{n}
\end{array}
\label{ass}
\end{equation}
$\Box .$
\vskip 0.3cm

This definition allows us to cast the bialgebra $H_\gamma (m,\Delta ,\eta
,\varepsilon )$ into a standard FRT bialgebra with deformed product. In this
way, higher tensor corepresentation of the bialgebra $H(m_\theta ,\Delta
,\eta ,\varepsilon ,R)$ are equivalent to non-commutative corepresentations
of $H_\gamma (m,\Delta ,\eta ,\varepsilon ).$

The twisting by 2-cocycles of quasitriangular Hopf algebras is due to
Drinfel'd \cite{drinf2} who has shown that starting from a quasitriangular
Hopf algebra, a new quasitriangular Hopf algebra is obtained twisting by a
2-cocycle the coproduct and the quasitriangular structure ${\cal R}$ . In
our case, we shall be interested in a partial twisting: we shall need just a
twisting of the product or a twisting of the dual-quasitriangular structure,
but not both together. We shall show that the braiding introduced by the
non-commutative coaction boils down to a partial twisting of the usual FRT
bialgebras, which in general do not preserves dual-quasitriangularity. To
this end, we extract some dual results from the Drinfeld analysis.

Following ref. \cite{maj}, we introduce a 2-cocycle on the bialgebra $%
Hom(H^{\otimes 2},k)$ as being an invertible element of $H^{\otimes 2},$ in
the sense of the product $*,$ satisfying the condition $\phi
_{23}*((I_{H}\otimes \Delta )\circ \phi )=\phi _{12}*((\Delta \otimes
I_{H})\circ \phi )\,\,\,.$

The maps $\phi :H^{\otimes 2}\rightarrow k$ being a 2-cocycle give rise to a
new bialgebra structure on $H_{o}$, namely $H_{\phi }(m_{\phi },\Delta ,\eta
,\varepsilon ),$ with a twisted product $m_{\phi }$%
\begin{eqnarray}
m_{\phi } &:&H_{o}^{\otimes 2}\rightarrow H_{o}  \label{mphi0} \\
m_{\phi } &=&\phi *m*\overline{\phi }  \nonumber
\end{eqnarray}
Moreover, if $\phi $ is a bialgebra bicharater, i.e.,

\begin{equation}
\begin{array}{c}
\phi (m\otimes I_{H})=\phi _{13}*\phi _{23} \\
\\
\phi (I_{H}\otimes m)=\phi _{13}*\phi _{12}
\end{array}
\label{bich}
\end{equation}
then the 2-cocycle condition leads to the Quantum Yang-Baxter equation

\begin{equation}
\phi _{12}*\phi _{13}*\phi _{23}=\phi _{23}*\phi _{13}*\phi _{12}
\label{phiyb}
\end{equation}

The above proposition provides the framework to interpret the bialgebra $%
H_\gamma $ as twisted one. In fact, let us assume the map $\theta ,$
introduced in eq. (\ref{hmap}), can be written as

\begin{equation}
\theta (T{_j^k)=\theta }_{jn}^{km}T_m^n=\rho _j^m\,T_m^n\,\overline{\rho }%
_n^k  \label{grho}
\end{equation}
i.e., the map $\theta $ admits the factorization

\begin{equation}
{\theta }_{ij}^{kl}=\rho _{i}^{l}\,\overline{\rho }_{j}^{k}  \label{titarho}
\end{equation}
Then, we may introduce the bialgebra bicharacter $\phi :H_{o}^{\otimes
2}\rightarrow k,$ inherited from the associativity assignment (\ref{ass}),
defined by the relations
\begin{equation}
\begin{array}{c}
\phi (T_{i}^{k}\otimes T_{j}^{l})=(\varepsilon \otimes \rho
)\,\,(T_{i}^{k}\otimes T_{j}^{l}) \\
\\
\rho (T_{i}^{k})=\rho _{i}^{k}\, \\
\\
\phi (T_{i}^{k}\otimes e)=\phi (e\otimes T_{i}^{k})=1 \\
\\
\phi (m\otimes I_{H})=\phi _{13}*\phi _{23} \\
\\
\phi (I_{H}\otimes m)=\phi _{13}*\phi _{12}
\end{array}
\label{phi}
\end{equation}
with $\rho :H_{o}\rightarrow k$ an invertible map, i.e., there exist $%
\overline{\rho }$ such that $\overline{\rho }*\rho =\rho *\overline{\rho }=$
$\varepsilon $, and $e$ is unit of the algebra $H_{o}$ . Then $\phi $ is a
2-cocycle, giving rise to the twisted product on $H_{o}:$
\begin{equation}
\begin{array}{c}
m_{\phi }(T_{i}^{k}\otimes T_{j}^{l})=T_{i}^{k}\,\,\widetilde{T}_{j}^{l} \\
\\
\widetilde{T}_{i}^{j}=\rho _{i}^{m}\,T_{m}^{n}\,\overline{\rho }_{n}^{j}
\end{array}
\label{mphi}
\end{equation}
Observe that the condition (\ref{gamarel2}) is trivially fulfilled. Then,
the bialgebra $H_{\gamma }(m,\Delta ,\eta ,\varepsilon ,R)$, eq.(\ref{hgam}%
), is isomorphic to a partial twisting by the 2-cocycle $\phi ,$(\ref{phi}),
of the standard FRT bialgebra $H(m,\Delta ,\eta ,\varepsilon ,R)$, (\ref{rtt}%
). The twisting can be performed on the product, thus obtaining the
bialgebra isomorphism $H_{\gamma }(m,\Delta ,\eta ,\varepsilon ,R)=H(m_{\phi
},\Delta ,\eta ,\varepsilon ,R),$ with the quadratic ideal ${\frak B}M-M%
{\frak B}$ being expressed as
\begin{equation}
m_{\phi }^{op}(T_{i}^{k}\otimes T_{j}^{l})=(R*m_{\phi }*\overline{R\,}%
)\,(T_{i}^{k}\otimes T_{j}^{l})  \label{rphittt}
\end{equation}
or, alternatively, one may leaves the product untwisted, but apply the
twisting onto the $R$ map,
\[
R^{\phi }=\overline{\phi }_{21}*R*\phi
\]
and, again, $H_{\gamma }(m,\Delta ,\eta ,\varepsilon ,R)=H(m,\Delta ,\eta
,\varepsilon ,R^{\phi })$ with the ideal, eq. (\ref{rtt}), expressed as
\begin{equation}
m^{op}(T_{i}^{k}\otimes T_{j}^{l})=(R^{\phi }*m*\overline{R}\,^{\phi
})(T_{i}^{k}\otimes T_{j}^{l})  \label{phirtt}
\end{equation}

In general, both schemes spoil out dual-quasitriangularity. However, in
order to leave open the connection with statistical systems where Boltzman
weight plays the role of the $R$ matrix, and the monodromy matrix are the $%
T_{i}^{k}$, it would be relevant if $R$ still provides a representation of
the bialgebra $H_{\gamma },$ as in fact it happens with dual-quasitriangular
bialgebras. This is achieved if $R^{\phi }$ to satisfy the relations

\begin{equation}
\begin{array}{c}
R^{\phi }\,(m\otimes I_{H})=R_{13}^{\phi }*R_{23}^{\phi } \\
\\
R^{\phi }\,(I_{H}\otimes m)=R_{13}^{\phi }*R_{12}^{\phi }
\end{array}
\label{phiq}
\end{equation}
so that $R^{\phi }$ is a solution of the Quantum Yang-Baxter equation. In
this way, those quantum bialgebras arising from non-commutative
corepresentations with factorizable $\theta $ -map, can be mapped into
standard FRT bialgebras by a twisting of the original $R.$

\section{Integrability}

The deep relation between Hopf Algebras and two dimensional physical systems
stems from the integrability condition. As we saw in the previous section,
for bialgebras arising from non-commutative corepresentations with a
factorizable $\theta $ -map it is possible to arrive to dual
quasitriangularity structure, then connection goes as usual .We shall see in
this section in which way a general bialgebra $H_\gamma $ may be associated
to some integrable systems. The main question is the introduction of the
spectral parameter, which is related to the coupling constant of the
physical system. In doing so, we proceed as in ref. \cite{drinf} by
regarding a collection of vector spaces $V(\lambda )={\bf C}^n$ for every $%
\lambda \in {\bf C,}$ with basis $b(\lambda )=\{e_i(\lambda ),$ $%
i=1,...,n\}. $

For each value of the spectral parameter $\lambda $ the coordinate functions
$T_{i}^{j}$ generates a bialgebra ${\sl H}_{o}(\lambda )$ with coproduct $%
\Delta T_{i}^{j}(\lambda )=T_{i}^{k}(\lambda )\otimes T_{k}^{j}(\lambda )$
and counit $\epsilon (T_{i}^{j}(\lambda ))=\delta _{i}^{j}.$ Furthermore,
the union of the ${\sl H}_{o}(\lambda )$ of these bialgebras for all values
of $\lambda ,$ i.e., ${\cal H}_{o}=\cup _{\lambda }{\sl H}_{o}(\lambda ),$
is also a bialgebra with the same coproduct and counit. Also, a ${\cal H}%
_{o} $-comodule structure on $V(\lambda )$ is obtained by the coaction $%
\delta _{V(\lambda )}\,e_{i}(\lambda )=T_{i}^{k}(\lambda )\otimes
e_{i}(\lambda ).$

Now, let us considers the map

\begin{equation}
\begin{array}{c}
{\frak B}(\lambda ,\mu ):V(\lambda )\otimes V(\mu )\longrightarrow V(\mu
)\otimes V(\lambda ) \\
\\
e_i(\lambda )\otimes e_j(\mu )\longrightarrow B_{ij}^{kl}(\lambda ,\mu
)\,e_k(\mu )\otimes e_l(\lambda )
\end{array}
\label{blambda}
\end{equation}
and the QLS defined by the quadratic algebra

\begin{equation}
A=%
{\displaystyle {\oplus _\lambda V^{\otimes }(\lambda ) \over R}}%
\label{lambdasp}
\end{equation}
where

\[
{\frak R}=\cup _{\lambda ,\mu }[1\otimes 1-B(\lambda ,\mu )]\,\,\,V(\lambda
)\otimes V(\mu )
\]
In order to obtain a structure of ${\cal H}_o$ comodule on $A$ , we define,
following the previous section, the map $\gamma $ as

\begin{equation}
\begin{array}{c}
\gamma :V(\lambda )\otimes {\cal H}_o\longrightarrow {\cal H}_o\otimes
V(\lambda ) \\
\\
\gamma (e_i(\lambda )\otimes T_j^k{(\mu ))=\gamma }_{ijn}^{klm}(\lambda ,\mu
)T_l^n{(\mu )}\otimes {e}_m(\lambda )
\end{array}
\label{lambdag}
\end{equation}
and the coaction on the tensor product space $V(\lambda )\otimes V(\mu )$

\begin{equation}
\delta _{V\otimes V}^{\gamma }=(m\otimes I_{V(\mu )})(I_{H_{\gamma }}\otimes
\gamma (\lambda ,\mu )\otimes I_{V{(\mu )}})(\delta _{V(\lambda )}\otimes
\delta _{V(\mu )})  \label{lambcvv}
\end{equation}
In analogy with {\it proposition 2}, this can be extended to a map,

\[
\delta ^\gamma :A\longrightarrow {\cal H}_o\otimes A
\]
supplying $A$ with a left ${\cal H}_o$-comodule structure$.$

The results analogous to those ones of {\it propositions 1-3} of the section
2 are still valid provided the replacements

\begin{equation}
\begin{array}{c}
{\gamma }_{ijn}^{klm}(\lambda ,\mu )=\delta _i^m\theta _{jn}^{kl}(\lambda
,\mu ) \\
\\
\theta _{ij}^{pl}(\lambda ,\mu )\theta _{pk}^{rs}(\lambda ,\mu )-\delta
_j^s\theta _{ik}^{rl}(\lambda ,\mu )=0 \\
\\
\theta _{jn}^{kn}(\lambda ,\mu )=\delta _j^k
\end{array}
\label{lambgam}
\end{equation}
and

\begin{equation}
M_{ij}^{kl}(\lambda ,\mu )=T_i^k{(\lambda )\,\,}\theta _{jn}^{km}(\lambda
,\mu )T_m^n{(\mu )}  \label{lambdam}
\end{equation}
So that, the ideal ${\frak B}M-M{\frak B}$ of {\it proposition 2} becomes in

\begin{equation}
B_{ij}^{mn}(\lambda ,\mu )\,\,M_{mn}^{kl}(\lambda ,\mu )-M_{ij}^{mn}(\mu
,\lambda )\,B_{mn}^{kl}(\lambda ,\mu )  \label{lambbmb}
\end{equation}
In order to study integrability, we analyze this equation. Eq.(\ref{lambbmb}%
). In general, a $\gamma $ map just satisfying the condition of the {\it %
prop.1, }does not lead to integrability making necessary to impose
additional conditions on it. In the following we study some options.

Our first ansatz is to require that

\begin{equation}
{\gamma }_{imj}^{mkl}(\lambda ,\mu )\equiv \delta _i^l\ \theta
_{mj}^{mk}=\delta _i^l\delta _j^k  \label{gamtraz}
\end{equation}
In this way, assuming $B$ to be invertible, one can multiply (\ref{lambbmb})
by the inverse of $B$ and then make the contraction of the free index of
both $B$ and $B^{-1},$ thus reaching the integrability condition

\begin{equation}
T(\lambda )T(\mu )=T(\mu )T(\lambda )  \label{integ}
\end{equation}
Here, $T(\lambda )$ means the trace $T_m^m(\lambda ).$

There is a less obvious way to recover integrability. Condition (\ref
{lambbmb}) in terms of $\theta $ is,
\begin{equation}
B_{ij}^{kl}(\lambda ,\mu )\ T_k^m(\lambda )\theta _{lv}^{nu}(\lambda ,\mu
)T_u^v(\mu )=\ T_i^k(\lambda )\theta _{jv}^{lu}(\lambda ,\mu )T_u^v(\mu
)B_{kl}^{mn}(\lambda ,\mu )  \label{bgttpar}
\end{equation}
If the following non trivial commutation holds, $,$%
\begin{equation}
B_{ij}^{rk}(\lambda ,\mu )\theta _{sr}^{sl}(\lambda ,\mu )=\theta
_{si}^{sr}(\lambda ,\mu )B_{rj}^{lk}(\lambda ,\mu )  \label{bggb}
\end{equation}
we may now contract (\ref{bgttpar}) with $\theta _{ar}^{ai}(\lambda ,\mu )$
and, after using (\ref{bggb}), we get

\begin{equation}
\begin{array}{l}
B_{rj}^{il}(\lambda ,\mu )\theta _{si}^{sk}(\lambda ,\mu )\
T_{k}^{m}(\lambda )\theta _{lv}^{nu}(\lambda ,\mu )T_{u}^{v}(\mu )= \\
\\
\,\,\,\,\,\,\,\,\,\,\,\,\,\,\,\,\,\,\,\,\,\,\,\,\,\,\,\,\,\,\,\,\,\,\,\,\,\,%
\,\,\,\,\,\,\,\,\,\,\,\,\,\,\,\,\,\,\,\,\,\,\,\,\,\,\,\,\,\,\,\,\,\,\,\theta
_{sr}^{si}(\lambda ,\mu )\ T_{i}^{k}(\lambda )\theta _{jv}^{lu}(\lambda ,\mu
)T_{u}^{v}(\mu )B_{kl}^{mn}(\lambda ,\mu )
\end{array}
\end{equation}
and now we proceed as in the previous case: multiplying by $%
(B^{-1})_{ab}^{rj}$ , and then performing the contractions $(a,m)\,,$ $%
(b,n), $ thus getting

\begin{equation}
\theta _{sa}^{sk}(\lambda ,\mu )\ T_{k}^{a}(\lambda )\theta
_{bv}^{bu}(\lambda ,\mu )T_{u}^{v}(\mu )=\theta _{sk}^{si}(\lambda ,\mu )\
T_{i}^{k}(\lambda )\theta _{lv}^{lu}(\lambda ,\mu )T_{u}^{v}(\mu )
\label{bgttpar2}
\end{equation}
Introducing the quantum matrix $\widetilde{T}$ of the final of the previous
section, the new integrability condition is written as

\begin{equation}
\widetilde{T}(\lambda )\widetilde{T}(\mu )=\widetilde{T}(\mu )\widetilde{T}%
(\lambda )  \label{integ2}
\end{equation}
This last approach to integrability may have an interpretation in the
framework of statistical models through their monodromy matrices $%
T_{i}^{k}(\lambda )$. In those models, periodic boundary conditions drive to
the transfer matrix by taking the trace over the auxiliary space of the
monodromy matrix, which means a sum over all the edge states (see for
example ref. \cite{marti}). The objects $\widetilde{T}(\lambda )$ means a
weighted sum over these edge states, so that $\theta $ seems to behave as a
twisting factor on the boundary conditions. The periodic ones correspond to
the trivial choice $\theta _{ij}^{kl}=\delta _{i}^{l}\delta _{j}^{k},$ and
we speculate that many other kind of boundary conditions would be reached by
a suitable choice of $\theta $ \cite{skly}$.$

In the next section we present a multiparametric example constructed from
the quantum plane and fulfilling the first integrability condition.

\section{The Quantum Plane}

Let us consider the quantum plane $A_{q}^{2\mid 0}$ described by

\begin{equation}
e_1e_2=q\ e_2e_1  \label{qpq}
\end{equation}

In the basis $\{e_{1}\otimes e_{1,}e_{1}\otimes e_{2},e_{2}\otimes
e_{1},e_{2}\otimes e_{1}\}$, this relation can be expressed by means the
quadratic form $B$ as

\begin{equation}
B=\left[
\begin{array}{cccc}
1 & 0 & 0 & 0 \\
0 & 0 & q & 0 \\
0 & q & 1-q^2 & 0 \\
0 & 0 & 0 & 1
\end{array}
\right]  \label{bqp}
\end{equation}
which is a solution of the Yang-Baxter equation

\begin{equation}
B_{12}B_{23}B_{12}=B_{23}B_{12}B_{23}  \label{yb}
\end{equation}

It is worth remarking that the following construction leads to the same
structure for other choice of $B$, as the symmetric and idempotent $%
B^{\prime },$

\begin{equation}
B^{\prime }=\frac 1{q+q^{-1}}\left[
\begin{array}{cccc}
q+q^{-1} & 0 & 0 & 0 \\
0 & q-q^{-1} & 2 & 0 \\
0 & 2 & q^{-1}-q & 0 \\
0 & 0 & 0 & q+q^{-1}
\end{array}
\right]  \label{bpq1}
\end{equation}

This $B^{\prime }$ is not a solution of Yang-Baxter equation but, in the
Manin construction for pseudo-symmetric quantum space\cite{man}, it enables
to characterize all the endomorphism of the quantum plane by the relation $%
B^{\prime }M-MB^{\prime }$ as the only solution to the master relation $%
(I-B^{\prime })M(I+B^{\prime }).$

The endomorphism matrix $T$ is

\begin{equation}
T=\left[
\begin{array}{cc}
a & b \\
c & d
\end{array}
\right]  \label{tqp}
\end{equation}
We find a multiparametric $\theta (r,p,s),$ solution of the coassociativity,
counity and the integrability condition , eqs. (\ref{gamarel2},\ref{gamarel3}%
,\ref{gamtraz}), that it is factorizable
\begin{equation}
\theta _{ij}^{kl}=\rho _{i}^{l}\,\overline{\rho }_{j}^{k}  \label{tpp}
\end{equation}
with
\begin{equation}
\rho =\left[
\begin{array}{ccc}
1 &  & r/s \\
&  &  \\
-s/p &  & (1-r)/p
\end{array}
\right]  \label{rho}
\end{equation}
and $\overline{\rho }$ is it inverse.

This means that the induced map $\phi =\varepsilon \otimes \rho $ is a
2-cocycle whenever $s\neq 0,$ hence there is an obstruction in to obtain the
undeformed limit $s\rightarrow 0.$ In this sense, this $\phi (r,p,s)$ is not
a 2-cocycle for the whole spectrum of its parameters and the twisting
arising from it does not yields a continuous deformation of the algebra $%
H_{o}.$

From this matrix we obtain the rules to commute $e$ and $T,$

\begin{equation}
\gamma (e_i\otimes T_j^k)=\theta _{jp}^{ks}\ T_s^p\otimes e_i  \label{te}
\end{equation}
As we see, there are now three new deformation parameters $(p,r,s)$ beside
to $q$, which was introduced by the quadratic algebra of the quantum plane (%
\ref{qpq}).

The relations $BM-MB$ can be written in a compact form by introducing the
objects $\widetilde{T}_i^j=\theta _{ir}^{js}\ T_s^r$, such that

\begin{equation}
\widetilde{T}=\left[
\begin{array}{cc}
\widetilde{a} & \widetilde{b} \\
\widetilde{c} & \widetilde{d}
\end{array}
\right]  \label{ts}
\end{equation}
so that we get

\begin{equation}
\begin{array}{c}
a\,\widetilde{c}-q\,c\,\widetilde{a}=0 \\
\\
a\,\widetilde{\,b}-q\,b\,\widetilde{a}=0 \\
\\
b\,\widetilde{c}-c\,\widetilde{b}=0 \\
\\
c\,\widetilde{d}-q\,d\,\widetilde{c}=0 \\
\\
b\,\widetilde{d}-q\,d\,\widetilde{b}=0 \\
\\
a\,\widetilde{d}-d\,\widetilde{a}+(q^{-1}-q)\,c\,\widetilde{b}=0
\end{array}
\label{tts}
\end{equation}
and

\begin{equation}
\begin{array}{c}
\gamma (e_{i}\otimes \widetilde{a})=\widetilde{a}\otimes e_{i} \\
\\
\gamma (e_{i}\otimes \widetilde{d})=\widetilde{d}\otimes e_{i} \\
\\
\gamma (e_{i}\otimes \widetilde{b})=p\,\widetilde{b}\otimes e_{i} \\
\\
\gamma (e_{i}\otimes \widetilde{c})=p^{-1}\,\widetilde{c}\otimes e_{i}
\end{array}
\label{emqp}
\end{equation}
The relations (\ref{tts}) acquires a highly non trivial form in terms of the
$T_{i}^{k}.$ In this way, one can define a four parameter deformation of the
$M(2),$ namely $M_{q,p,r,s}(2),$

\begin{equation}
M_{q,p,r,s}(2)=\frac{k[T_{i}^{j}]}{I(BM-MB)}  \label{mq2}
\end{equation}
Also, the Grassmannian plane $A_{q}^{0\mid 2}(\xi _{1},\xi _{2})$, defined
by the relation
\[
\xi _{1}\xi _{2}=-%
{\displaystyle {1 \over q}}%
\xi _{2}\xi _{1}
\]
is naturally a $M_{q,p,r,s}(2)$-comodule \cite{man1}. This allows us to
define a determinant for this $M_{q,p,r,s}(2)$ from the coaction $\delta
_{V\otimes V}^{\gamma }$ on the object $\xi _{1}\xi _{2}$

\begin{equation}
\delta _{V\otimes V}^\gamma (\xi _1\xi _2)=D\otimes \xi _1\xi _2
\label{coacd}
\end{equation}
Thus we get

\begin{equation}
D=M_{12}^{12}-qM_{12}^{21}=a\,\widetilde{d}-q\,b\,\widetilde{c}=a\,d-\frac{q%
}{p}(1-r)\ b\,c-%
{\displaystyle {r \over s}}%
\ a\,c-q%
{\displaystyle {s \over p}}%
\ b\,d  \label{det}
\end{equation}
As was explained above, this is a non-perturbative deformation: the limit to
the undeformed case can't be taken simultaneously. However, there is a
sequential limit leading to another 3 and 2 parameters deformation. In fact
, if we take first the limit $r\rightarrow 0$ we get $M_{q,p,s}(2)$ and now
the remaining $\phi (p,s)$ is now a genuine 2-cocycle, so that the twisting
is well defined in the whole spectrum of $p$ and $s.$ Taking now a second
undeformed limit, we set $s\rightarrow 0,$ then we recover the biparametric $%
M_{q,p}(2)$ obtained by Manin {\it et al,\cite{man2}, }as non-standard
quantum groups. In these limits $\theta $ becomes in

\begin{equation}
\theta (p)=\left[
\begin{array}{cccc}
1 & 0 & 0 & 0 \\
0 & 0 & p & 0 \\
0 &
{\displaystyle {1 \over p}}%
& 0 & 0 \\
0 & 0 & 0 & 1
\end{array}
\right]  \label{tita2}
\end{equation}
and now, the 2-cocycle $\phi =\epsilon \otimes \rho $ becomes in

\begin{equation}
\rho =\left[
\begin{array}{cc}
1 & 0 \\
0 & 1/p
\end{array}
\right]  \label{rho1}
\end{equation}
and the relations $BM-MB,$ or $R^{\phi }*m=m^{op}*R^{\phi },$ reduce to

\begin{equation}
\begin{array}{r}
a\,c-pq\,c\,a=0 \\
\\
a\,b-p^{-1}q\,b\,a=0 \\
\\
b\,c-p^2c\,b=0 \\
\\
c\,d-p^{-1}q\,d\,c=0 \\
\\
b\,d-pq\,d\,b=0 \\
\\
a\,d-d\,a+p(q^{-1}-q)\,c\,b=0
\end{array}
\label{bmqp}
\end{equation}
which define a two parametric $M_{q,p}(2)$ as the quotient algebra

\begin{equation}
M_{q,p}(2)=\frac{k[T_{i}^{j}]}{I(BM-MB)}  \label{mq3}
\end{equation}
for $i$ and $j$ from 1 to 2. It is worth remarking that $R^{\phi }(q,p)=%
\overline{\phi }_{21}(p)*R(q)*\phi (p)$ is a solution of the Quantum Yang
Baxter equation, and assuming $R^{\phi }$ is a bialgebra bicharacter, it
supplies $M_{q,p}(2)$ with a dual-quasitriangular structure.

This $\gamma $ give rise to the following relations between matrix entries
and the coordinates of the quantum plane,

\begin{equation}
\begin{array}{c}
\gamma (e_{i}\otimes a)=a\otimes e_{i} \\
\\
\gamma (e_{i}\otimes d)=d\otimes e_{i} \\
\\
\gamma (e_{i}\otimes b)=pb\otimes e_{i} \\
\\
\gamma (e_{i}\otimes c)=p^{-1}c\otimes e_{i}
\end{array}
\label{emqp1}
\end{equation}
Defining $\widetilde{T}_{i}^{j}=\theta _{il}^{jk}\ T_{k}^{l}$ we get

\begin{equation}
\widetilde{T}=\left[
\begin{array}{ccc}
a &  & pc \\
&  &  \\
{\displaystyle {b \over p}}%
&  & d
\end{array}
\right]  \label{tsomqp}
\end{equation}
and with the coproduct

\begin{equation}
\triangle \widetilde{T}_i^j=\widetilde{T}_i^k\otimes \widetilde{T}_k^j
\label{cop-coact}
\end{equation}
The determinant becomes in

\[
D=\det\nolimits_{q,p}=ad-p^{-1}qbc
\]
that satisfy the following commutation relations

\[
\begin{array}{r}
Da-aD=0 \\
\\
Db-p^{-2}bD=0 \\
\\
Dc-p^2cD=0 \\
\\
Dd-dD=0
\end{array}
\]
With these properties, and assuming that $D$ is an invertible element of $%
M_{q,p}(2)$, the antipode can be defined

\begin{equation}
S\left[
\begin{array}{cc}
a & b \\
c & d
\end{array}
\right] =D^{-1}\left[
\begin{array}{cc}
d & -(pq)^{-1}b \\
-pqc & a
\end{array}
\right]  \label{antip}
\end{equation}
and, consequently, the Hopf algebra $GL_{q,p}(2)$ is obtained . Thus, this
biparametric deformation of $M(2)$ developed in ref. \cite{man2} can be
alternatively interpreted in the framework of non-commutative
corepresentations.

\section{Concluding Remarks}

We have introduced a new ingredient in the theory of corepresentations of
Quantum (semi)Groups admitting non-commutativity between endomorphism matrix
entries and quantum space coordinates. This feature give rise to an extra
deformation of all the involved structures. Our approach is not a full
braiding as those obtained from the quasitensor category of $k$-modules\cite
{maj}. In a less ambitious project, we have just redefined the coaction
showing that it is possible to introduce a non trivial map $\gamma :V\otimes
H_{\gamma }\longrightarrow H_{\gamma }\otimes V$ without spoiling out the
Hopf algebra and the comodule structures provided that $\gamma $ turns $%
H_{\gamma }$ into a Quantum Matrix Group. No additional modification were
introduced in the usual structure of bialgebras, preserving the product and
coproduct untouched. However, provided a factorizable $\gamma ,$ this
generalization of the coaction boils down to a twisting of the algebra
structure or a twist of the $R$-matrix, and in some cases it is posible to
recover a quasitriangular FRT bialgebra.

Although the map $\gamma $ seems too constrained, we have found non trivial
solutions introducing many new deformation parameters, still under the
additional condition of integrability. This last point was also analyzed,
showing that integrability can be reached at least in two independent ways.
Working on the Quantum Planes $A_{q}^{2\mid 0}$and $A_{q}^{0\mid 2}$ as
examples, it was shown that the biparamatric deformation of $GL(2,{\bf C),}$
namely $GL_{pq}(2),$ can be regarded as coacting in a twisted way over the
standard quantum plane.

\section{Acknowledgments}

We are greatly indebted to M. L. Bruschi for enlightening discussions. Also,
the authors thank CONICET-Argentina for financial support. R. T. also thanks
Universit\'{e} Mediterrane\'{e} (Aix-Marseille II) for finantial support.


\end{document}